\documentclass{article}

\begin{document}

\title{A short proof of the error term in Simpson's rule}
\markright{A short proof of the error term in Simpson's rule}
\author{Hajrudin Fejzi\'{c}}

\maketitle

When the integral $\int_a^b f(x)\,dx$ is difficult to evaluate by the
standard integration techniques, we can use numerical integration
formulas such as Simpson's rule to approximate the integral.  It is
well-known that the error in Simpson's rule is related to the fourth
derivative of the integrand.  In addition to proofs using techniques
from numerical analysis, the typical proofs using basic ideas from
calculus involve repeated applications of the Mean Value Theorem (see
\cite{G} or \cite{O}) or Taylor's Theorem (see \cite{T}, which also
includes a detailed discussion of this error estimate.)  We offer here
a modified version of this proof.

As indicated in \cite{G} and \cite{T}, it is sufficient to focus on
the special case $n=2$ and then use the intermediate property of
derivatives for the general formula.  Suppose that $f$ is continuous
on $[a,b]$ and four times differentiable on $(a,b)$.  Simpson's rule
is obtained by integrating the quadratic polynomial $q$ that agrees
with the function $f$ at the points $a$, $b$, and $c=(a+b)/2$.  An
elementary computation reveals that
$$
\int_a^b \bigl(f(x)-q(x)\bigr)\,dx =
\int_a^bf(x)\,dx - \frac{b-a}{6}\bigl(f(a)+4f(c)+f(b)\bigl).
$$
Let $E$ be the value of the above expression; our goal is to show that
$$
E = -\frac{f^{iv}(\xi)}{32\cdot90}(b-a)^5 =
\frac{f^{iv}(\xi)}{4!}\int_a^b xp(x)\,dx,
$$
where $p(x)=(x-a)(x-c)(x-b)$ and $\xi$ is some point
between $a$ and $b$.  (We have omitted the elementary computation that
gives the second equality.)

Let $k$ be a real number and consider the continuous function
$$
\phi(x) = f(x)-q(x)-kp(x)-\frac{E}{\int_{a}^{b} xp(x)\,dx}\,xp(x).
$$
Since $f$ and $q$ agree at the points $a$, $b$, and $c$, the function
$\phi$ is zero at these points.  A simple computation shows that
$\int_a^b p(x)\,dx=0$ and it follows that $\int_a^b
\phi(x)\,dx=0$. Since $\int_a^c p(x)\,dx\neq 0$, we can select $k$ so
that $\int_a^c \phi(x)\,dx=0$.  From the equation
$$
0 = \int_a^b \phi(x)\,dx = \int_a^c \phi(x)\,dx +
\int_c^b \phi(x)\,dx,
$$
it follows that $\int_c^b \phi(x)\,dx=0$.  Since the function $\phi$
is continuous, there exist points $u\in(a,c)$ and $v\in(c,b)$ such
that $\phi(u)=0$ and $\phi(v)=0$.  Hence, the function $\phi$ is zero
at five distinct points in the interval $[a,b]$.  Repeated application
of Rolle's Theorem yields a point $\xi$ in $(a,b)$ such that
$\phi^{iv}(\xi)=0$ and thus
$$
f^{iv}(\xi)- \frac{4!E}{\int_{a}^{b} xp(x)\,dx} = 0,
$$
completing the proof.

\end{document}